\input amstex
\input Amstex-document.sty

\pageno 323

\def\shorttitle{Collapsed Riemannian Manifolds with Bounded Sectional Curvature}
\def\shortauthor{Xiaochun Rong}

\topmatter
\title\nofrills{\boldHuge Collapsed Riemannian Manifolds with Bounded Sectional
Curvature*}
\endtitle

\author \Large Xiaochun Rong$^\dag$\endauthor

\thanks *Supported partially by NSF Grant DMS 0203164 and a research found from
Beijing Normal University. \endthanks

\thanks $^\dag$Rutgers University, New Brunswick, NJ 08903, USA and Beijing Normal
University, Beijing 100875, China. E-mail: rong\@math.rutgers.edu
\endthanks

\abstract\nofrills \centerline{\boldnormal Abstract}

\vskip 4.5mm

{\ninepoint In the last two decades, one of the most important
developments in Riemannian geometry is the collapsing theory of
Cheeger-Fukaya-Gromov. A Riemannian manifold is called (sufficiently)
collapsed if its dimension looks smaller than its actual dimension while
its sectional curvature remains bounded (say a very thin flat torus
looks like a circle in a bared eyes).

We will survey the development of collapsing theory and its applications
to Riemannian geometry since 1990. The common starting point for all of
these is the existence of a singular fibration structure on collapsed
manifolds. However, new techniques have been introduced and tools from
related fields have been brought in. As a consequence, light has been
shed on some classical problems and conjectures whose  statements do not
involve collapsing. Specifically, substantial progress has been made on
manifolds with nonpositive curvature, on positively pinched manifolds,
collapsed manifolds with an a priori diameter bound, and subclasses of
manifolds whose members satisfy additional topological conditions e.g.
$2$-connectedness.

\vskip 4.5mm

\noindent {\bf 2000 Mathematics Subject Classification:} 53C.

\noindent {\bf Keywords and Phrases:} Gromov-Hausdorff convergence,
Collapsing, F-structure.}
\endabstract
\endtopmatter

\document

One of the most important developments in Riemannian geometry over
the last two decades is the structure theory of Cheeger-Fukaya-Gromov
for manifolds $M^n$ of bounded sectional curvature, say
$|\text{sec}_{M^n}|\leq 1$,
which are sufficiently collapsed. Roughly, $M^n$ is called
$\epsilon$-collapsed, if it appears to have dimension less than $n$,
unless the metric is rescaled by a factor $\geq\epsilon^{-1}$.

 For scaling reasons, collapsing and boundedness of
tend to oppose one another.
Nevertheless, very collapsed manifolds with bounded curvature do in fact exist.
 For example, a very thin cylinder is very collapsed, although its
curvature vanishes identically.

If one fixes $\epsilon$ and in addition, a bound, $d$, on the
diameter, then in each dimension, there only finitely many manifolds,
which are not $\epsilon$-collapsed; see [Ch]. The basic result of collapsing
theory states the existence of a constant $\epsilon(n)>0$, such
that a manifold which is $\epsilon$-collapsed, for $\epsilon\leq
\epsilon(n)$, has a particular kind of singular fibration structure
with flat (or ``almost flat'') fibers. The fibers lie in the
$\epsilon$-collapsed directions; see [CG1,2], [CFG], [Fu1-3].

The first nontrivial example of a collapsing sequence
 with bounded curvature (described in more detail below)
 was constructed by M. Berger in 1962; see [CFG].  The
first major result on the collapsed manifolds (still a cornerstone
of the theory) was M. Gromov's characterization of ``almost flat
manifolds'' i.e. manifolds admitting a sequence of metrics with
curvature and diameter going to zero. Gromov showed that such manifolds
are infranil. Later in [Ru], they were shown to actually
be nilmanifolds; compare [GMR].

We will survey the development of collapsing theory and its
applications to Riemannian geometry since 1990; compare [Fu4].
The common starting point for all of these is the above
mentioned singular fibration structure. However,
new techniques have been introduced and
tools from related fields have been brought in.
As a consequence, light has been shed on some classical problems and
conjectures whose  statements do not involve collapsing. Specifically,
substantial progress has been made on manifolds with nonpositive
curvature, on positively pinched manifolds, collapsed manifolds
with an a priori diameter bound, and subclasses of manifolds whose members
satisfy additional topological conditions e.g. $2$-connectedness.

\vskip2mm

\specialhead \noindent \boldLARGE 1. Collapsed manifolds of
bounded sectional curvature
\endspecialhead

\vskip2mm

Convention: unless otherwise specified, ``collapsing'' refers to
a sequence of Riemannian manifolds with sectional
curvature bounded in absolute value by one and injectivity
radii uniformly converge to zero, while  ``convergence'' means ``convergence
with respect to the Gromov-Hausdorff distance.

Recall that a map
from a metric space $(X,d_X)$ to a metric space $(Y,d_Y)$ is called
an \text{$\epsilon$-Gromov-Hausdorff} approximation, if
of $f(X)$ is $\epsilon$-dense in $Y$ and
if $|d_X(x,x')-d_Y(f(x),f(x'))|<\epsilon$. The
Gromov-Hausdorff distance between two (compact) metric spaces is the
infimum of $\epsilon$ as above,
 for all possible $\epsilon$-Gromov-Hausdorff
approximations from $X$ to $Y$ and vice versa.
(To be more precise, one should say ``pseudo-distance'', since
isometric metric spaces have distance zero.)
The collection of all compact metric spaces
is complete with respect to the Gromov-Hausdorff distance.

\vskip2mm

{\noindent \bf a. Flat manifolds, collapsing by scaling and torus
actions}

\vskip2mm

 For fixed $(M,g)$, the family, $\{(M,\epsilon^2g)\}$
converges to a point as $\epsilon\to 0$. However, if the
curvature is not identically zero,
it blows up. On the other hand,
 for any compact flat manifold, $(M,g)$, the
the manifolds, $(M,\epsilon^2g)$ continue to be flat.
 More generally, if $(M,g)$ is a (possibly nonflat) manifold with an
isometric torus $T^k$-action for which all $T^k$-orbits have
the same dimension, then one obtains a collapsing sequence by rescaling
$g$ along the orbits i.e. by putting
 $g_\epsilon=\epsilon^2g_0\oplus g_0^\perp$,
where $g_0$ is the restriction of $g$ to the tangent space of
a $T^k$-orbit and $g_0^\perp$ is the orthogonal complement.
A computation shows that $g_\epsilon$ has bounded sectional curvature
independent of $\epsilon$. The collapse constructed by Berger in 1962
was of this type. In his example, $M^3$ is the unit
$3$-sphere and the $S^1$ action is by rotation in the fibers of the
Hopf fibration $S^1\to S^3\to S^2$. The limit of this collapse
is the 2-sphere with a metric of constant curvature $\equiv 4$; see
[Pet].

More generally,
a collapsing construction has been given
by Cheeger-Gromov for manifolds which admit
certain mutually compatible local torus actions (possibly
by tori of different dimensions) for which all orbits
have positive dimension; see the notion of
{\it F-structure} given below and (1.2.1). As above,
for each individual local torus
action, one obtains locally defined collapsing sequence.
The problem is to patch together these local collapsings. If the orbits
are not all of the same dimension, the patching
requires a suitable scaling of the metric
(by a large constant) in the transition regions between
orbits of different dimensions; see [CG1]. Hence, in contrast to
the Berger example, in general the diameters of such nontrivially
patched collapsings necessarily go to infinity.

\vskip2mm

{\noindent \bf b. Almost flat manifolds and collapsing by
inhomogeneous scaling}

\vskip2mm

Although a compact nilmanifold (based on a nonablian nilpotent Lie group)
admits no flat metric,
a sequence metrics on such a manifold which collapses to a point
can be constructed by a suitable inhomogeneous scaling process;
see [Gr1]. As an example, regard a compact nilmanifold
$M^3$ as the total space of a principle circle bundle over
a torus. A canonical
metric $g$ on $M^3$ splits into horizontal and vertical complements,
$g=g_h\oplus g_h^\perp$.
Then $g_\epsilon=(\epsilon g_h)\oplus (\epsilon^2 g_h^\perp)$ has
bounded sectional curvature independent of $\epsilon$,
while $(M^3,g_\epsilon)$ converges to a point. The inhomogeneity
of the scaling is essential in order for the curvature to remain bounded;
compare Theorem 3.4.

\vskip2mm

{\noindent\bf c. Positive rank F-structure and collapsed
manifolds}

\vskip2mm

The notion of an F-structure may be viewed as a generalization of that of
a torus action. An F-structure $\Cal F$ on a manifold is defined by
an atlas $\Cal F=\{(V_i,U_i,T^{k_i})\}$, satisfying the following
conditions:

\noindent (1.1.1) $\{U_i\}$ is a locally finite open cover for $M$.

\noindent (1.1.2) $\pi_i: V_i\to U_i$ is a finite normal covering and
$V_i$ admits an effective torus $T^{k_i}$-action such that it extends
to a $\pi_1(U_i)\ltimes T^{k_i}$-action.

\noindent (1.1.3) If $U_i\cap U_j\ne \emptyset$, then $\pi_i^{-1}(U_i\cap U_j)$
and $\pi_j^{-1}(U_i\cap U_j)$ have a common finite covering on which
the lifting $T^{k_i}-$ and $T^{k_j}$-actions commute.

If $k_i=k$, for all $i$, then $\Cal F$ is called {\it pure}. Otherwise,
$\Cal F$ is called {\it mixed}. The compatibility condition, (1.1.3),
implies that $M$
decomposes into {\it orbits}. (an orbit at a point is the smallest set
containing all the projections of the $T^{k_i}$-orbits at the point.)
The minimal dimension of all such
orbits is called the {\it rank of $\Cal F$}.
An orbit is called {\it regular}, if it has a tubular neighborhood
in which the orbits form a fibration. Otherwise, it
is called {\it singular}. An F-structure $\Cal F$ is called {\it
polarized} if all $T^{k_i}$-actions are almost free. An F-structure is
called {\it injective} (resp. semi-injective) if the inclusion of any
orbit to $M$ induces an injective (resp. nontrivial) map on the
fundamental groups.

A {\it Cr-structure} is an injective F-structure
with an atlas that
satisfies two additional properties: i) $V_i=D_i\times T^{k_i}$ and
$T^{k_i}$ acts on $V_i$ by the multiplication.
ii) If $U_i\cap U_j\ne \emptyset$,
then $k_i<k_j$ or vice versa; see [Bu1]. This notion arises in
the context of nonpositive curvature.

A metric is called an {\it $\Cal F$-invariant} (or simply invariant),
if the local $T^{k_i}$-actions are isomeric. For any F-structure, there
exists an invariant metric.

A manifold may not admit
any nontrivial F-structure; compare Corollary 2.5.  In fact, a simple
necessary condition for a closed manifold $M^{2n}$ to admit a positive rank
F-structure is the vanishing of its Euler characteristic; see [CG1].

A necessary and sufficient condition for the existence of a collapsing
sequence of metrics is the existence of an F-structure of positive rank;
see [CG1], [CG2].

\vskip2mm

\proclaim{Theorem 1.2 (Collapsing and F-structure of positive
rank)}{\rm  ([CG1,2])} Let $M$ be a manifold without boundary.

\noindent (1.2.1) If $M$ admits a positive rank (resp. polarized)
F-structure, then $M$ admits a continuous one-parameter family of
invariant metrics $g_\epsilon$ such that $|\text{sec}_{g_\epsilon}|
\le 1$ and the injectivity radius (resp. volume) of $g_\epsilon$
converges uniformly to zero as $\epsilon\to 0$.

\noindent (1.2.2) There exists a constant $\epsilon(n)$ (the critical
injectivity radius) such that if $M^n$ admits a metric $g$ with
$|\text{sec}_g|\le 1$ and the injectivity radius is less than
$\epsilon(n)$ everywhere, then $M$ admits a positive rank F-structure
almost compatible with the metric.
\endproclaim

\vskip1mm

The F-structure in (1.2.2) is actually a substructure of
a so called {\it nilpotent Killing structure on $M$} whose orbits are
infra-nilmanifolds; see [CFG] and compare to Theorem 3.5. Such an infra- nilmanifold
orbit at a point contains {\it all} sufficiently collapsed directions
of the metric;
the orbit of its sub F-structure, which is defined by the `center'
of the infra-nilmanifold, only contains the most collapsed directions
comparable to the injectivity radius
at a point. A unsolved problem pertaining to nilpotent structures is
whether a collapse as in (1.2.1) can be constructed
for which the diameters of the nil-orbits
converge uniformly to zero (as holds for F-structures).

The construction of the F-structure in (1.2.2) relies only
on the local geometry. Hence, (1.2.2) can be applied to a collapsed
region in a complete manifold of bounded
sectional curvature. In this way, for such a manifold, one obtains
a {\it thick-thin} decomposition, in which the thin
part carries an F-structure of positive rank; see [CFG].

Theorem 1.2
has been the starting point for many subsequent investigations
of collapsing in various situations. The guiding principle is
that additional geometrical properties of a collapsing should
be mirrored in properties of its associated F-structure,
which in turn, puts constraints on the topology.
For instance, if a collapsing satisfies additional geometrical conditions
such as: i) volume small, ii) uniformly bounded diameter, iii) nonpositive
curvature, iv) positive pinched curvature, v) {\it bounded covering geometry}
i.e. the injectivity radii of the Riemannian universal covering has a uniform
positive lower bound,
then one may expect corresponding
 topological properties of the F-structure such as: i) existence of a
polarization,
ii) pureness, iii) existence of a Cr-structure, iv) the existence of a circle orbit, v) injective
F-structure. Results on such correspondences and their applications
will occupy the rest of this paper.

\vskip2mm

{\noindent \bf d.  Topological invariants associated to a volume
collapse}
 \vskip2mm

The existence of a sufficiently (injectivity radius) collapsed metric as in (1.2.2) imposes
constraints on the underlying topology. For instance, the
{\it simplicial volume} of
$M$ vanishes; see [Gr3]. As mentioned earlier, for a closed $M^{2n}$, the
Euler characteristic of $M^{2n}$ also vanishes; see [CFG].

In this subsection, we focus on some topological invariants associated to
certain (partially) volume collapsed metrics: the {\it minimal volume},
the {\it $L^2$-signature} and the {\it limiting $\eta$-invariant};
see below.

The minimal volume,  $\text{MinVol}(M)$, of $M$,
 is the infimum of the volumes over all
complete metrics with $|\text{sec}_M|\le 1$. Clearly, $\text{MinVol}(M)$
is a topological invariant.
Gromov conjectured that there exists a constant $\epsilon(n)>0$ such
that $\text{MinVol}(M^n)<\epsilon(n)$ implies that $\text{MinVol}(M^n)=0$
(the gap conjecture for minimal volume).
By Theorem 1.2, it would suffice to show that a sufficiently volume
collapsed manifold admits a polarized F-structure. On a $3$-manifold,
any positive rank F-structure has a polarized substructure and thus
Theorem 1.2 implies Gromov's gap conjecture in dimension 3. However,
for $n\ge 4$, there are $n$-manifolds which admit a positive rank F-structure
but which admit no polarized F-structure; see [CG1].

\vskip2mm

\proclaim{Theorem 1.3 (Volume collapse and Polarized F-structure)}
{\rm ([Ro2])} There is a constant $\epsilon>0$ such that if
$\text{MinVol}(M^4)<\epsilon$, then $M^4$ admits a polarized
F-structure and thus $\text{MinVol}(M^4)=0$.
\endproclaim

\vskip1mm

 For a complete open manifold with
 bounded sectional curvature and finite volume
(necessarily volume collapsed near infinity), the integral of an invariant
polynomial of the curvature form may depend on the particular metric;
see [CG3]. It is of interest
to find a class of metrics for which integral of characteristic forms
have a topological interpretation.
Cheeger-Gromov showed that for any open
complete manifold $M^{4k}$ of finite volume
and bounded covering geometry outside some compact subset, the
integral of the Hirzebruch signature form over $M^{4k}$ is independent of the
metric; see [CG3] and the references therein. Cheeger-Gromov showed that
this integral is equal to the so called  $L_2$-signature and
conjectured that it can take only rational values.
(The notion of $L_2$-signature, whose definition involves the concept of
Von Neumann dimension, was first introduced by Atiyah and Singer in
the context of coverings of compact manifolds.)

\vskip2mm

\proclaim{Theorem 1.4 (Rationality of geometric signature)} {\rm
([Ro3])} If an open complete manifold, $M^4$,  of finite volume
has bounded covering geometry outside a compact subset, then the
integral of the Hirzebruch signature form over $M^4$ is a rational
number.
\endproclaim

\vskip1mm

The main idea is to show that $M^4$ admits a polarized F-structure
$\Cal F$ outside some compact subset and an
exhaustion by compact submanifolds, $M_i^4$, such that the restriction of
$\Cal F$ to the boundary of $M^n_i$ is injective. The integral over $M^4$
is the limit of the
integrals over $M^4_i$, to which we apply the Atiyah-Patodi-Singer formula
to reduce to showing the rationality of the limit of the $\eta$-invariant terms.
By making use of the special property of $\Cal F$ and Theorem 1.5 below,
we are able to conclude that the limit of the $\eta$-invariant term is rational.

Cheeger-Gromov showed that if a sequence of volume collapsed metrics
on a closed manifold $N^{4n-1}$ have bounded covering geometry,
then the sequence of the associated $\eta$-invariants converges and
the limit is independent of the particular sequence of such metrics. They
conjectured that the limit is rational.

\vskip2mm

\proclaim{Theorem 1.5 (Rationality of limiting $\eta$-invariants)}
{\rm ([Ro1])} If a closed manifold $N^3$ admits a sequence of
volume collapsed metrics with bounded covering geometry, then
$N^3$ admits an injective F-structure and the limit of the
$\eta$-invariants is rational.
\endproclaim

\vskip1mm

The idea is to show that $N^3$ admits an injective F-structure
$\Cal F$. For an injective F-structure, the collapsing constructed
in (1.2.1) has bounded covering geometry and may be used to
compute the limit. Results from $3$-manifold topology play
a role in the proof of the existence of the injective F-structure.

\vskip2mm

\specialhead \noindent \boldLARGE 2. Collapsed manifolds with
nonpositive sectional curvature
\endspecialhead

\vskip2mm

A classical result of Preismann says that for a closed manifold $M^n$ with
negative sectional curvature, any abelian subgroup of the fundamental group
is cyclic. By
bringing in the discrete group technique, Margulis showed that if the metric
is normalized such that
$-1\le \text{sec}_{M^n}\le 0$, then there exists at least one point at which
the injectivity radius is bounded below by a constant $\epsilon(n)>0$.

The study of the subsequent study of collapsed manifolds with
$-1\le \text{sec}\le 0$ may be viewed as
an attempt to describe the special
circumstances under which the conclusions of
the Preismann and Margulis theorem can fail, if the hypothesis is weakened
to nonpositive curvature; see [Bu1-3], [CCR1,2], [Eb], [GW], [LY], [Sc].

 A collapsed metric with nonpositive curvature tends
to be rigid in a precise sense; see (2.2.1) and (2.2.2).
Namely, there exists a {\it canonical} Cr-structure
whose orbits are flat totally geodesic submanifolds.  Of necessity,
the construction
of this Cr-structure is global.  By contrast, the construction of
less precise  (but more generally existing) F-structure is local; see [CG2].

Let $M^n=\tilde M^n/\Gamma$, where $\tilde M^n$ denotes the universal
covering space of $M^n$ with the pull-back metric.
A {\it local splitting structure}
on a Riemannian manifold is a
$\Gamma$-equivariant assignment to each point (of an open dense subset of
$\tilde M^n$) a specified neighborhood and a specified
isometric splitting of this neighborhood, with a nontrivial
Euclidean factor. Hence, a necessary condition for a local splitting
structure is the existence of a plane of zero curvature,
at every point of $M^n$. A local splitting structure is {\it abelian}
if the projection to $M^n$ of every nontrivial Euclidean factor as above
is a closed {\it embedded} flat submanifold, and n addition,
if two projected leaves
intersect, then one of them is contained in the other.

\vskip2mm

\proclaim{Theorem 2.1 (Abelian local splitting structure and
Cr-structure)} {\rm ([CCR1])} Let $M^n$ be a closed manifold of
$-1\le \text{sec}_{M^n} \le 0$.

\noindent (2.1.1) If the injectivity radius is smaller than $\epsilon(n)>0$ everywhere, then $M^n$ admits an abelian local splitting structure.

\noindent (2.1.2) If $M^n$ admits an abelian local splitting structure, then it admits a compatible Cr-structure, whose orbits are the
flat submanifolds (projected leaves) of the abelian local splitting structure.
In particular, $\text{MinVol}(M^n)=0$.
\endproclaim

\vskip1mm

Theorem 2.1 was conjectured by Buyalo, who proved the cases $n=3,
4$; see [Bu1--3], [Sc].

Let $\tilde x\in M^n$.
 Let $\Gamma_\epsilon(\tilde x)\ne 1$ denote
the subgroup of $\Gamma$ generated by those $\gamma$
whose displacement function, $\delta_\gamma( \tilde x)=
d(x,\gamma(\tilde x))$, satisfies
$d(x,\gamma(\tilde x))<\epsilon$. (In the application, $\epsilon$ is small.)
If all
$\Gamma_\epsilon(\tilde x)$ are abelian, then the minimal sets, $\{\text{Min}
(\Gamma_\epsilon(\tilde x))\}$,
of the $\Gamma_\epsilon(\tilde x)$ give the desired
abelian local splitting structure in (2.1.1).
In general,
$\Gamma_\epsilon(\tilde x)$ is only Bieberbach. Then, a crucial ingredient in
(2.1.1) is the existence of a `canonical' abelian subgroup
of $\Gamma_\epsilon(\tilde x)$
of finite index consisting of those elements which are
{\it stable} in the sense of [BGS].
In spirit, the proof
of (2.1.2) is similar to the construction in [CG2], but the techniques
used are quite
different.

The following are some specific questions pertaining to abelian local splitting structures:

\noindent (2.2.1) If some metric $g$ on $M$ of nonpositive sectional
curvature has an abelian local splitting structure, does
every nonpositively curved metric also have such a structure?
\vskip2mm

\noindent (2.2.2) If $M$ has a Cr-structure, does every any nonpositively
curved metric on $M$ have a compatible local splitting structure?
\vskip2mm

Note that an affirmative answer to (2.2.1) and (2.2.2)
would imply a kind of {\it semirigidity}.
It would imply that all nonpositively curved metrics on
$M$ are alike in a precise sense.

\vskip2mm

\proclaim{Theorem 2.3 (F-structure and local splitting structure)} {\rm ([CCR2])} Let $X^n$, $M^n$ be closed
manifolds such that $X^n$ admits a nontrivial $F$-structure.  Let $f:X^n\to M^n$ have nonzero degree. Then every
metric of nonpositive sectional curvature on $M^n$ has a local splitting structure.
\endproclaim

\vskip1mm

We conjecture that if an F-structure has positive rank, then the
local splitting structure is abelian. This conjecture, whose proof
would provide an affirmative answer to
(2.2.2), has been verified in dimension 3 and in some additional
 special cases; see [CCR2].

We conclude this section with two consequences of Theorem 2.3.

\vskip2mm

\proclaim{Corollary 2.4 (Generalized Margulis Lemma)} {\rm
([CCR2])} Let $M^n$ be a closed manifold of nonpositive sectional
curvature. If the Ricci curvature is negative at some point, then
for every metric with $|\text{sec}|\le 1$, there is a point with
injectivity radius $\ge \delta(n)>0$.
\endproclaim

\vskip1mm

Another consequence is a geometric obstruction for
a nontrivial F-structure.

\vskip2mm

\proclaim{Corollary 2.5 (Nonexistence of F-structure)} {\rm
([CCR2])} If a closed manifold $M$ admits a metric of nonpositive
sectional curvature such that the Ricci curvature is negative at
some point, then $M$ does not admit a nontrivial F-structure.
\endproclaim

\vskip2mm

\specialhead \noindent \boldLARGE 3. Collapsed manifolds with
bounded sectional curvature and diameter
\endspecialhead

\vskip2mm

In this section, we discuss the class of collapsed manifolds of bounded
sectional curvature whose diameters are also bounded. By the Gromov's
compactness theorem, any sequence of such collapsed manifolds
contains a convergent subsequence; see [GLP]. Hence, without loss of the
generality, we only need to consider convergent collapsing sequences.

\noindent (3.1) Let $M^n_i @>d_{GH}>> X$ denote a sequence of closed
manifolds converging to a {\it compact} metric space $X$ such
that $|\text{sec}_{M^n_i}|\le 1$ and $\dim(X)<n$.

\vskip2mm

\example{Main Problem 3.2} For $i$ large,
investigate relations between geometry and topology of
$M^n_i$ and that of $X$. The following are some specific problems
and questions.

\noindent (3.2.1) Find topological obstructions for the existence of
$M^n_i$ as in (3.1).

\noindent (3.2.2) To what extent is the topology of the
 $M^n_i$ in (3.1) stable when $i$ is sufficiently large?

\noindent (3.2.3) Under what additional conditions is it
true that $\{M^n_i\}$ as in (3.1)
contains a subsequence of constant diffeomorphism type? If all $M^n_i$ are
diffeomorphic, then to what extent do the metrics converge?
\endexample

\vskip1mm

Note that by the Cheeger-Gromov convergence theorem,
the above problems are well understood
in the noncollapsed situation $\dim(X)=n$.

\vskip2mm

\proclaim{Theorem 3.3 (Convergence)} {\rm ([Ch], [GLP])} Let
$M^n_i@>d_{GH}>> X$ be as in (3.1) except $\dim(X)=n$. Then for
$i$ large, $M^n_i$ is diffeomorphic to some fixed $M^n$ which is
homeomorphic to $X$ and there are diffeomorphisms, $f_i: M^n\to
M^n_i$, such that the pulled back metrics,  $f_i^*(g_i)$, converge
to a metric, $g_\infty$, in the $C^{1,\alpha}$-topology
$(0<\alpha<1$).
\endproclaim

\vskip1mm

Note that as a consequence of Theorem 3.3, topological stability
of a sequence as in (3.1)
will immediately yield a corresponding finiteness result in terms
of the dimension and bounds on curvature and diameter.

\vskip2mm

{\noindent \bf e. Structure of collapsed manifolds with bounded
diameter}

\vskip2mm

As described in Section 1, any
closed nilmanifold $M^n$ admits metrics collapsing to a point.

\vskip2mm

\proclaim{Theorem 3.4 (Almost flat manifolds)} {\rm ([Gr1])} Let
$M^n_i@>d_{GH}>> X$ be as in (3.1). If $X$ is a point, then a
finite normal covering space of $M^n_i$ of order at most $c(n)$ is
diffeomorphic to a nilmanifold $N^n/\Gamma_i$ ($i$ large), where
$N^n$ is the simply connected nilpotent group.
\endproclaim

\vskip1mm

Theorem 3.4 can be promoted to a description of
convergent collapsing sequence, of manifolds, $M^n_i$,
as in (3.1). As mentioned following
Theorem 1.2, any sufficiently collapsed
manifold admits a nilpotent Killing structure; see [CFG].
Here a bound on diameter forces the nilpotent Killing structure
to be pure.

For a closed Riemannian manifold $M^n$, its frame bundle $F(M^n)$ admits a
canonical metric determined by the Riemannian connection up to
a choice of a bi-invariant metric on $O(n)$. A fibration,
$N/\Gamma\to F(M^n)\to Y$, is called $O(n)$-invariant if the $O(n)$-action
on $F(M^n)$ preserves both the fiber $N/\Gamma$ (a nilmanifold) and
the structural group. By the $O(n)$-invariance, $O(n)$ also acts on
the base space $Y$. A canonical metric is invariant if
its restriction on each $N/\Gamma$ is left-invariant.
A {\it pure nilpotent Killing structure}
on $M$ is an $O(n)$-invariant fibration
on $F(M^n)$ for which the canonical metric is also invariant.

\vskip2mm

\proclaim{Theorem 3.5 (Fibration)} {\rm ([CFG])} Let
$M^n_i@>d_{GH}>> X$ be as in (3.1). Then $F(M^n_i)$ equipped with
canonical metrics contains a convergent subsequence,
$F(M^n_i)@>d_{GH}>>Y$, and $F(M^n_i)$ admits an $O(n)$-invariant
fibration $N/\Gamma_i\to F(M^n_i) \to Y$ for which the canonical
metric is $\epsilon_i$-close in the $C^1$ sense to some invariant
metric, where $\epsilon_i\to 0$.
\endproclaim

\vskip1mm

The following properties are crucial for the study of particular instances of collapsing as in (3.1).

\vskip2mm

\proclaim{Proposition 3.6} Let $M^n_i@>d_{GH}>> X$ be as in (3.1).

\noindent (3.6.1) (Regularity) ([Ro5]) For any $\epsilon>0$, $M^n_i$
admits an invariant metric $g_i$ such that $\min (\text{sec}_{M^n_i})-\epsilon
\le \text{sec}_{(M^n_i,g_i)}\le \max (\text{sec}_{M^n_i})+\epsilon$ for $i$ large.

\noindent (3.6.2) (Equivariance) ([PT], [GK]) The induced $O(n)$actions on $Y$ from
the $O(n)$-action on $F(M^n_i)$ are $C^1$-close and therefore are all
$O(n)$-equivariant for $i$ large.
\endproclaim

\vskip2mm

{\noindent \bf f. Obstructions to collapsing with bounded
diameter}

\vskip2mm

\proclaim{Theorem 3.7 (Polarized F-structure and vanishing minimal
volume)}\linebreak {\rm ([CR2])} Let $M^n_i@>d_{GH}>> X$ be as in
(3.1). Then the F-substructure associated to the pure nilpotent
Killing structure on $M^n_i$ contains a (mixed) polarized
F-structure. In particular, $\text{MinVol}(M^n_i)=0$.
\endproclaim

\vskip1mm

Theorem 3.7 may be viewed as a weak version of the Gromov's gap conjecture.
Note that the associated F-structure on $M^n_i$  may not be
polarized. The existence of a polarized substructure puts
constraints on the singularities of the structure.

\vskip2mm

\proclaim{Theorem 3.8 (Absence of symplectic structure)} {\rm
([FR3])} Let $M^n_i@>d_{GH}>> X$ be as in (3.1). If $\pi_1(M^n_i)$
is finite, then $M^n_i$ does not support any symplectic structure.
\endproclaim

\vskip1mm

The proof of Theorem 3.8 includes a nontrivial extension of the
well known fact that any $S^1$-action on a closed simply connected symplectic manifold
which preserves the symplectic structure has a nonempty fixed
point set.

A geometric obstruction to the existence of a collapsing sequence in (3.1) is
provided by:

\vskip2mm

\proclaim{Theorem 3.9 (Geometric collapsing obstruction)} {\rm
([Ro7])} Let $M^n_i@>d_{GH}>> X$ be as in (3.1). Then $\limsup
(\max_{M^n_i} \text{Ric}_{M^n_i})\ge 0$.
\endproclaim

\vskip1mm

A key ingredient in the proof is a generalization
of a theorem of Bochner asserting that
a closed manifold of negative Ricci curvature admits no nontrivial
invariant pure F-structure (Bochner's original
theorem only guarantees the nonexistence of a nontrivial
isometric torus action.)

\vskip2mm

\proclaim{Theorem 3.10 (Pure injective F-structure)} {\rm ([CR1])}
Let $M^n_i@>d_{GH}>> X$ be as in (3.1). If $M^n_i$  has bounded
covering geometry and $\pi_1(M^n_i)$ is torsion free, then for $i$
large $M^n_i$ admits a pure injective F-structure.
\endproclaim

\vskip1mm

{\noindent \bf g. The topological and geometric stability}

\vskip2mm

In this subsection, we address Problems (3.2.2) and (3.2.3).
Observe that by the Gromov's Betti number estimate, [Gr2], the
sequence in (3.1)
contains a subsequence whose cohomology {\it groups},
$H_*(M^n_i,\Bbb Q)$, are all isomorphic. On the other hand,
examples have been found showing that $\{H^*(M^n_i,\Bbb Q)\}$ can
contain infinitely many distinct {\it ring} structures; see [FR2].

\vskip2mm

\proclaim{Theorem 3.11 ($\pi_q$-Stability)} {\rm ([FR2]; compare
[Ro4], [Tu])} Let $M^n_i @ > d_{GH} >> X$ be as in (3.1). Then for
$q\ge 2$ and after passing to a subsequence, the $q$-th homotopy
group $\pi_q(M^n_i)$ are all isomorphic, provided that
$\pi_q(M^n_i)$ are finitely generated (e.g. $\text{sec}_{M^n_i}\ge
0$ or $\pi_1(M^n_i)$ is finite).
\endproclaim

\vskip1mm

Note that in contrast to the Betti number bound, Theorem 3.11 does not hold
if upper bound on the sectional curvature is removed; see [GZ].

We now discuss sufficient topological conditions for diffeomorphism
stability. Consider the sequence of fibrations, $N/\Gamma_i\to
F(M^n_i)\to Y$, associated to (3.1). One would like to know
when all $N/\Gamma_i$ are diffeomorphic.

\vskip2mm

\proclaim{Proposition 3.12} {\rm ([FR4])} Let $M^n_i @ > d_{GH} >>
X$ be as in (3.1). If $\pi_1(M^n_i)$ contains no free abelian
group of rank two, then $N/\Gamma_i$ is diffeomorphic to a torus.
\endproclaim

\vskip1mm

In low dimensions, we have:

\vskip2mm

\proclaim{Theorem 3.13 (Diffeomorphism stability---low
dimensions)} {\rm ([FR3], [Tu])} For $n\le 6$, let $M^n_i @ >
d_{GH}
>> X$ be as in (3.1). If $\pi_1(M^n_i)=1$, then there is a
subsequence all whose members are diffeomorphic.
\endproclaim

\vskip1mm

Note that for $n\ge 7$, one cannot expect Theorem 3.13; see [AW]. Hence,
additional restrictions are required in higher dimensions.
 Observe that if $M^n_i$
are 2-connected, then all $T^k\to F(M^n_i)
\to Y$ are equivalent as principle $T^k$-bundles. In particular
all $F(M^n_i)$ are diffeomorphic.

Using (3.6.2), Petrunin-Tuschmann
showed that the equivalence can be chosen that is also
$O(n)$-equivariant, and concluded the diffeomorphism stability
for two-connected manifolds; see [PT]. For the special case in which
the $M^n_i$ are positively pinched, the same
conclusion was obtained independently in [FR1] via a different
approach.

We introduce a topological condition which when $M^n_i$ is simply
connected, reduces to the assumption that $\pi_2(M^n_i)$ is finite.
In the nonsimply connected case however, there are manifolds with
$\pi_2(M)$ infinite, which satisfy our condition.

 Let $\tilde M$ denote the universal
covering of $M$. For a homomorphism, $\rho: \pi_1(M)\to \text{Aut}(\Bbb Z^k)$,
the semi-direct product, $\tilde M\times_{\pi _1(M)} \Bbb Z^k$,
is a bundle of $\rho(\pi_1(M))$-modules which can be viewed as
 a {\it local coefficient system} over $M$.  We denote
 it
 by $\Bbb Z^k_\rho$. Let
$\tilde b_q(M,\Bbb Z^k_\rho)$
denote the rank of the cohomology group, $H^q(M,\Bbb Z^k_\rho)$,
with the local coefficient system $\Bbb Z^k_\rho$.
We refer to the integer
$$\tilde b_q(M,\Bbb Z^k)=\underset{\rho: \pi _1(M)\to \text{Aut}
(\Bbb Z^k) }\to \max\{\tilde b_q(M,\Bbb Z^k_\rho)\}$$
as the $q$-th {\it twisted Betti number of $M$}. Clearly,
$\tilde b_q(M,\Bbb Z^k)$
is a topological invariant of $M$.  Moreover,
$k\cdot b_2(M,\Bbb Z)\le
\tilde b_q(M,\Bbb Z^k)$, with equality if $\pi_1(M)=1$.
\vskip2mm

\proclaim{Theorem  3.14 (Diffeomorphism stability and geometric
stability)} {\rm ([FR4])} Let $M^n_i @ > d_{GH} >> X$ be as in
(3.1) with $k=n-\dim(X)$. Assume that $M^n_i$ satisfies:

\noindent (3.14.1) $\pi_1(M^n_i)$ is a torsion group with torsion
exponents uniformly bounded from above.

\noindent (3.14.2) The second twisted Betti number $\tilde b_2(M^n_i,
\Bbb Z^k)=0$.

\noindent Then there are diffeomorphisms, $f_i$,
 from $M^n$ to (a subsequence of)
$\{M^n_i\}$, such
that the distance functions of pullback metrics, $f^*_i(g_i)$,
on $M^n$, converge to a pseudo-metric $d_\infty$ in $C^0$-norm. Moreover,
$M^n$ admits a foliation with leaves diffeomorphic to
flat manifolds (that are not necessarily compact) and a vector
$V$ tangent to a leaf if and only if $\|V\|_{g_i}
\to 0$.
\endproclaim

\vskip1mm

The proof of Theorem 3.14 is quite involved.

Finally, we mention that J. Lott has systematically investigated
the analytic aspects for a collapsing in (3.1); for details,
see [Lo1-3].

\vskip2mm

\specialhead \noindent \boldLARGE 4. Positively pinched manifolds
\endspecialhead

\vskip2mm

In this section, we further investigate a subclass of the
class of collapsed manifolds with
bounded diameter: collapsed manifolds with pinched positive sectional
curvature; see [AW], [Ba], [Es], [P\"u] for examples.

In the spirit of Theorem 3.4, we first give the following classification
result.

\vskip2mm

\proclaim{Theorem 4.1 (Maximal collapse with pinched positive
curvature)}{\rm ([Ro8])} Let $M^n_i @ > d_{GH} >> X$ be as in
(3.1) such that $\text{sec}_{M^n_i}\ge \delta >0$. Then
$\dim(X)\ge \frac{n-1}2$ and equality implies that $\tilde
M^n_i\overset{\text{diffeo}}\to \simeq S^n/\Bbb Z_{q_i}$ (a lens
space), where $\tilde M^n_i\to M^n_i$ is a covering space of order
$\le\frac{n+1}2$.
\endproclaim

\vskip1mm

By Theorem 3.5, (3.6.1) and Proposition 3.12, the proof of Theorem 4.1
reduces to
the classification of   positively curved manifolds which admit
invariant pure F-structures of maximal rank; see [GS].

\vskip2mm

\proclaim{Theorem 4.2 (Positive pinching and almost cyclicity of
$\pi_1$)} {\rm ([Ro6])} Let $M^n_i @ > d_{GH} >> X$ be as in (3.1)
such that $\text{sec}_{M^n_i}\ge \delta>0$. Then for $i$
sufficiently large, $\pi_1(M^n_i)$ has a cyclic subgroup whose
index is less than $w(n)$.
\endproclaim

\vskip1mm

By Theorem 3.5 and (3.6.1), the following result easily implies Theorem 4.2.

\vskip2mm

\proclaim{Theorem 4.3 (Symmetry and almost cyclicity of $\pi_1$)}
{\rm ([Ro6])} Let $M^n$ be a closed manifold of positive sectional
curvature. If $M^n$ admits an invariant pure F-structure, then
$\pi_1(M^n)$ has a cyclic subgroup whose index is less than a
constant $w(n)$.
\endproclaim

\vskip1mm

In the special case of a free isometric action, from the
homotopy exact sequence associated to the
fibration, $S^1\to M^n\to M^n/S^1$, together with the
Synge theorem, one sees that
$\pi_1(M^n)$ is cyclic. The proof of the general case is
by induction on $n$ and is rather complicated.

We now consider the injectivity radius estimate. Klingenberg-Sakai
and Yau conjectured that the infimum of the injectivity radii of all
$\delta$-pinched metrics on $M^n$ is a positive number which
depends only on $\delta$ and the homotopy type of the manifold. By a result
of Klingenberg, this conjecture is easy in even dimensions. In odd
dimensions it is open.

\vskip2mm

\proclaim{Theorem 4.4 (Noncollapsing)} {\rm ([FR4]; compare [FR1],
[PT])} For $n$ odd, let $M^n$ be a closed manifold satisfying
 $0<\delta\le \text{sec}_{M^n}\le 1$ and
$|\pi_1(M^n)|\le c$. If $\tilde b(M^n,\Bbb Z^{\frac{n-1}2})=0$, then
the injectivity radius of $M^n$ is at least $\epsilon
(n,\delta,c)>0$.
\endproclaim

\vskip1mm

If Theorem 4.4 were false, then by Theorem  3.14 and (3.6.1) one could
assume the existence of
 a sequence, $(M,g_i)@ > d_{GH} >> X$, with $\delta/2\le
\text{sec}_{g_i}\le 1$,  such that the distance functions of
the metrics $g_i$
also converge. In view of the following theorem this would lead
to a contradiction.

\vskip2mm

\proclaim{Theorem 4.5 (Gluing)} {\rm ([PRT])} Let $(M,g_i)@ >
d_{GH}
>> X$ as in (1.3). If the distance functions of $g_i$ converge to
a pseudo-metric, then $\liminf (\min\text{sec}_{g_i})\le 0$.
\endproclaim

\vskip1mm

Let $f_i: (M,g_i)@>d_{GH}>>X$ denote an $\epsilon_i$
Gromov-Hausdorff approximation, where
$\epsilon_i\to 0$.
 For an open cover $\{B_j\}$ for $X$ by small (contractible) balls, the
assumption on the distance functions implies (roughly) that the tube,
$C_{ij}=f_i^{-1}(B_j)$, is a subset of $M$ independent of $i$. Clearly,
the universal covering $\tilde C_{ij}$ of $C_{ij}$ is noncompact. The idea is
to glue together the limits of the
 $\tilde C_{ij}$ (modulo some suitable
group of isometries with respect to the pullback metrics)
to form a noncompact metric space with curvature bounded
below by $\liminf (\min\text{sec}{}_{g_i})$ in the comparison sense; see [BGP], [Pe].
On the other hand,
the positivity of the curvature implies that the space so obtained would
have to be compact.

The above results on $\delta$-pinched manifolds may shed a light on the
topology of positively curved manifolds. It is
tempting to make the following conjecture (which seems very difficult).

\vskip2mm

\example{Conjecture 4.6} Let $M^n$ denote a closed manifold of
positive sectional curvature.

\noindent (4.6.1) (Almost cyclicity) $\pi_1(M^n)$ has a cyclic subgroup
with index bounded by a constant depending only on $n$.

\noindent (4.6.2) (Homotopy group finiteness) For $q\ge 2$,
$\pi_q(M^n)$ has only finitely many possible isomorphism classes
depending only on $n$ and $q$.

\noindent (4.6.3) (Diffeomorphism finiteness) If $\pi_q(M^n)=0$
($q=1, 2$), then $M^n$ can have only finitely many possible diffeomorphism types
depending only on $n$.
\endexample

\vskip1mm

Note that (4.6.1)--(4.6.3) are false for nonnegatively curved spaces.
By the results in this section, Conjecture 4.6 would follow from an
affirmative answer to the following:

\vskip2mm

\example{Problem 4.7 (Universal pinching constant)} ([Be], [Ro5]) Is there a constant
$0<\delta(n)<<1$ such that any closed $n$-manifold of positive
sectional curvature admits a $\delta(n)$-pinched metric?
\endexample

\vskip2mm

A partial verification of (4.6.2) is obtained by [FR2].
\vskip2mm

\proclaim{Theorem 4.8} {\rm ([FR2])} Let $M^n$ denote a closed
manifold of positive sectional curvature. For $q\ge 2$, the
minimal number of generators for $\pi_q(M^n)$ is less than
$c(q,n)$.
\endproclaim

\vskip1mm

Previously, by Gromov the minimal number of generators of
$\pi_1(M^n)$ is bounded above by a constant depending only on $n$.

\specialhead \noindent \boldLARGE References \endspecialhead

\widestnumber\key{APSSS}

\ref \key AW \by S. Aloff; N. R. Wallach \pages 93--97 \paper An
infinite family of $7$-manifolds admitting positive curved
Riemannian structures \jour Bull. Amer. Math. Soc. \vol 81 \yr
1975
\endref

\ref
\key BGS
\by W. Ballmann; M. Gromov; Schroeder
\pages
\paper Manifolds of nonpositive curvature
\jour Basel:
\newline Birkh\"auser, Boston, Basel, Stuttgart,
\vol
\yr 1985
\endref

\ref
\key Ba
\by Ya. V. Baza\u{\i}kin Y
\pages 1068--1085
\paper On a family of $13$-dimensional closed Riemannian manifolds of
positive curvature
\jour Sibirsk. Mat. Zh. 37 (in Russian), ii; English
translation in Siberian Math. J.
\vol 6
\yr 1996
\endref

\ref
\key Be
\by M. Berger
\pages
\paper Riemannian geometry during the second half of the twentieth century
\jour University lecture series
\vol 17
\yr 2000
\endref

\ref \key BGP \by Y. Burago; M. Gromov; Perel'man \pages 3--51
\paper A.D. Alexandov spaces with curvature bounded below \jour
Uspekhi Mat. Nauk \yr 1992 \vol 47:2
\endref

\ref \key Bu1 \by  S. Buyalo \pages 1135--1155 \paper Collapsing
manifolds of nonpositive curvature I \jour Leningrad Math. J.,
\vol 5 \yr 1990
\endref

\ref \key Bu2 \by  S. Buyalo \pages 1371--1399 \paper Collapsing
manifolds of nonpositive curvature II \jour Leningrad Math. J.,
\vol 6 \yr 1990
\endref

\ref
\key Bu3
\by  S. Buyalo
\pages 1--26
\paper Three dimensional manifolds with $Cr$-structure,
\jour Some Questions of Geometry in the Large, A.M.S. Translations
\vol 176
\yr 1996
\endref

\ref \key CCR1 \by J. Cao; J. Cheeger; X. Rong \pages 139--167
\paper Splittings and Cr-structure for manifolds with nonpositive
sectional curvature \jour Invent. Math. \yr 2001 \vol 144
\endref

\ref
\key CCR2
\by J. Cao; J. Cheeger; X. Rong
\pages
\paper Partial rigidity of nonpositively curved manifolds
\jour
\yr To appear
\vol
\endref

\ref \key Ch \by J. Cheeger \pages 61--75 \paper Finiteness
theorems for Riemannian manifolds \jour Amer. J. Math. \yr 1970
\vol 92
\endref

\ref \key CFG \by J. Cheeger; K. Fukaya; M. Gromov \pages 327--372
\paper Nilpotent structures and invariant metrics on collapsed
manifolds \jour J. A.M.S \yr 1992 \vol 5
\endref

\ref \key CG1 \by J. Cheeger; M. Gromov \pages 309--364 \paper
Collapsing Riemannian manifolds while keeping their curvature
bound I \jour J. Diff. Geom \yr 1986 \vol 23
\endref

\ref
\key CG2
\by  J. Cheeger; M. Gromov
\pages 269--298
\paper Collapsing Riemannian manifold while keeping their
curvature bounded II
\jour J. Differential Geom
\yr 1990
\vol 32
\endref

\ref \key CG3 \by J. Cheeger; M. Gromov \pages 115--154 \paper On
the characteristic numbers of complete manifolds of bounded
curvature and finite volume \jour H. E. Rauch Mem Vol I (Chavel
and Farkas, Eds) Springer, Berlin \yr 1985
\endref

\ref \key CG4 \by J. Cheeger; M. Gromov \pages 1--34 \paper Bounds
on the von Neumann dimension of $L_2$-cohomology and the
Gauss-Bonnet theorem for open manifolds \jour J. Diff. Geom \yr
1985 \vol 21
\endref

\ref \key CR1 \by J. Cheeger; X. Rong \pages 141--163 \paper
Collapsed manifolds with bounded diameter and bounded covering
geometry \jour Geome.Funct. Anal \yr 1995 \vol 5 No. 2
\endref

\ref \key CR2 \by  J. Cheeger; X. Rong \pages 411--429 \paper
Existence of polarized F-structure on collapsed manifolds with
bounded curvature and diameter \jour Geome. Funct. Anal \yr 1996
\vol 6, No.3
\endref

\ref
\key Eb
\by P. Eberlein
\pages 23--29
\paper A canonical form for compact nonpositively curved manifolds
whose fundamental groups have nontrivial center
\jour Math. Ann.
\vol 260
\yr 1982
\endref

\ref \key Es \by J.-H Eschenburg \pages 469--480 \paper New
examples of manifolds with strictly positive curvature \jour
Invent. Math \yr 1982 \vol 66
\endref

\ref \key FR1 \by  F. Fang; X. Rong \pages 641--674 \paper
Positive pinching, volume and homotopy groups \jour Geom. Funct.
Anal \yr 1999 \vol 9
\endref

\ref \key FR2 \by F. Fang; X. Rong \pages 135--158 \paper
Curvature, diameter, homotopy groups and cohomology rings \jour
Duke Math. J. \yr 2001 \vol 107 No.1
\endref

\ref \key FR3 \by  F. Fang; X. Rong \pages 75--86 \paper Fixed
point free circle actions and finiteness theorems \jour Comm.
Contemp. Math \yr 2000 \vol
\endref

\ref
\key FR4
\by  F. Fang; X. Rong
\pages
\paper The twisted second Betti number and convergence of collapsing
Riemannian manifolds
\jour To appear in Invent. Math
\yr 2002
\vol
\endref

\ref \key Fu1 \by K. Fukaya \pages 139--156 \paper Collapsing
Riemannian manifolds to ones of lower dimension \jour J. Diff.
Geome \yr 1987 \vol 25
\endref

\ref \key Fu2 \by K. Fukaya \pages 333--356 \paper Collapsing
Riemannian manifolds to ones of lower dimension II \jour J. Math.
Soc. Japan \yr 1989 \vol 41
\endref

\ref \key Fu3 \by K. Fukaya \pages 1--21 \paper A boundary of the
set of the Riemannian manifolds with bounded curvature and
diameters \jour J. Diff. Geome \yr 1988 \vol 28
\endref

\ref
\key Fu4
\by K. Fukaya
\paper Hausdorff convergence of Riemannian manifolds and its applications
\jour Recent Topics in Differential and Analytic Geometry (T. Ochiai, ed),
Kinokuniya, Tokyo
\yr 1990
\endref

\ref \key GMR \by P. Granaat; M. Min-Oo; E. Ruh \pages 1305--1312
\paper Local structure of Riemannian manifolds \jour Indiana Univ.
Math. J \yr 1990 \vol 39
\endref

\ref
\key GW
\by D. Gromoll; J. Wolf
\pages 545--552
\paper Some relations between the metric structure and the algebraic
structure of the fundamental group in manifolds of nonpositive
curvature
\jour Bull. Am. Math. Soc.,
\vol 4
\yr 1977
\endref

\ref \key Gr1 \by M. Gromov \pages 231--241 \paper Almost flat
manifolds \jour J. Diff. Geom \yr 1978 \vol 13
\endref

\ref \key Gr2 \by M. Gromov \pages 179--195 \paper Curvature
diameter and Betti numbers \jour Comment. Math. Helv \yr 1981 \vol
56
\endref

\ref \key Gr3 \by M. Gromov \pages 213--307 \paper Volume and
bounded cohomology \jour I.H.E.S. Pul. Math. \yr 1983 \vol 56
\endref

\ref
\key GLP
\by M. Gromov, J. Lafontaine; P.Pansu
\paper Structures metriques pour les varietes riemannienes
\jour CedicFernand Paris
\yr 1981
\endref

\ref \key GK \by K. Grove; H. Karcher \pages 11--20 \paper How to
conjugate $C^1$-close actions \jour Math. Z \vol 132 \yr 1973
\endref

\ref \key GS \by K. Grove, C. Searle \pages 137--142 \paper
Positively curved manifolds with maximal symmetry-rank \jour J.
Pure Appl. Alg \yr 1994 \vol 91
\endref

\ref \key GZ \by K. Grove; W. Ziller \pages 331--367 \paper
Curvature and symmetry of Milnor spheres \jour Ann. of Math \vol
152 \yr 2000
\endref

\ref
\key LY
\by B. Lawson; S. T. Yau
\pages
\paper On compact manifolds of nonpositive curvature
\jour J. Diff. Geom.,
\vol 7
\yr 1972
\endref

\ref
\key Lo1
\by J. Lott
\pages
\paper Collapsing and differential form Laplacian: the case
of a smooth limit space
\jour Duke Math. J
\vol
\yr To appear
\endref

\ref
\key Lo2
\by J. Lott
\pages
\paper Collapsing and differential form Laplacian: the case
of a smooth limit space
\jour  Preprint
\vol
\yr 2002
\endref

\ref
\key Lo3
\by J. Lott
\pages
\paper Collapsing and Dirac-type operators
\jour Geometriae Delicata, Special issue on partial differential equations and their
applications to geometry and physics
\vol
\yr To appear
\endref

\ref
\key Pe
\by G. Perel'man
\pages
\paper A.D. Alexandrov spaces with curvature bounded below II
\jour preprint
\vol
\yr
\endref

\ref
\key Pet
\by P. Petersen
\pages
\paper Riemannian geometry
\jour GTM, Springer-Verlag Berlin Heidelberg New York
\vol 171
\yr 1997
\endref

\ref \key PRT \by A. Petrunin; X. Rong; W. Tuschmann \pages
699--735 \paper Collapsing vs. positive pinching \jour Geom.
Funct. Anal
 \yr 1999
\vol 9
\endref

\ref
\key PT
\by A. Petrunin; W. Tuschmann
\pages
\paper Diffeomorphism finiteness, positive pinching, and second
homotopy
\jour Geom. Funct. Anal
\yr 1999
\vol 9
\endref

\ref \key P\"u \by T. P\"uttmann \pages 631--684 \paper Optimal
pinching constants of odd dimensional homogeneous spaces \jour
Invent. Math \vol 138 \yr 1999
\endref

\ref \key Ro1 \by  X. Rong \pages 535--568 \paper  The limiting
eta invariant of collapsed $3$-manifolds \jour J. Diff. Geom \yr
1993 \vol 37
\endref

\ref \key Ro2 \by  X. Rong \pages 475--502 \paper The existence of
polarized F-structures on volume collapsed $4$-manifolds \jour
Geom. Funct. Anal \yr 1993 \vol 3, No.5
\endref

\ref \key Ro3 \by  X. Rong \pages 513--554 \paper Rationality of
geometric signatures of complete $4$-manifolds \jour Invent. Math
\yr 1995 \vol 120
\endref

\ref \key Ro4 \by  X. Rong \pages 427--435 \paper Bounding
homotopy and homology groups by curvature and diameter \jour Duke.
Math. J \yr 1995 \vol 2
\endref

\ref \key Ro5 \by  X. Rong \pages 397--411 \paper On the
fundamental group of manifolds of positive sectional curvature
\jour Ann. of Math \yr 1996 \vol 143
\endref

\ref \key Ro6 \by  X. Rong \pages 47--64 \paper The almost
cyclicity of the fundamental groups of positively curved manifolds
\jour Invent. Math \yr 1996 \vol 126
\endref

\ref \key Ro7 \by  X. Rong \pages 381--392 \paper A Bochner
Theorem and applications \jour Duke Math. J \yr 1998 \vol 91, No.2
\endref

\ref \key Ro8 \by  X. Rong \pages 335--358 \paper Collapsed
manifolds with pinched positive sectional curvature \jour J. Diff.
Geom \yr 1999 \vol 51 No.2
\endref

\ref \key Ru \by E. Ruh \pages 1--14 \paper  Almost flat manifolds
\jour J. Diff. Geome. \yr 1982 \vol 17
\endref

\ref
\key Sc
\by V. Schroeder
\pages 19--26
\paper Rigidity of nonpositively curved graph manifolds
\jour Math. Ann.
\vol 274
\yr 1986
\endref

\ref \key Tu \by  W. Tuschmann \pages 413--420 \paper Geometric
diffeomorphism finiteness in low dimensions and homotopy
finiteness \jour Math. Ann \yr 2002 \vol 322
\endref

\enddocument